\documentclass[12pt]{amsart}
\usepackage{graphicx}
\usepackage{color}
\usepackage{setspace}
\usepackage{eufrak}
\usepackage{float}
\usepackage{etex}
\usepackage[all]{xy}
\usepackage{amssymb, amsmath}
\usepackage{ascmac}

\def\qed{\hfill $\Box$}
\def\proof{\noindent {\sl Proof} :\;  }
\def\t{\noindent}

\newcommand{\V}{\mathcal{V}}

\newcommand{\A}{\mathcal{A}}

\newcommand{\Proj}{\mathbb{P}}

\newcommand{\R}{\mathbb{R}}

\def\qed{\hfill $\Box$}
\def\proof{\noindent {\sl Proof} :\;  }
\def\t{\noindent}
\def\rd{\partial}

\def\bu{\mbox{\boldmath $u$}}

\def\b0{\mbox{\boldmath $0$}}

\def\I{{\rm I}}
\def\II{{\rm I\hspace{-.1em}I}}
\def\III{{\rm I\hspace{-.1em}I\hspace{-.1em}I}}
\def\IV{{\rm I\hspace{-.1em}V}}
\def\V{{\rm V}}
\def\VI{{\rm VI}}

\newtheorem{thm}{\bf Theorem}[section]

\newtheorem{prop}[thm]{\bf Proposition}

\newtheorem{rem}[thm]{\bf Remark}

\begin{document}
\title[Asymptotic Binary Differential Equations on surfaces]
{Binary differential equations 
at parabolic and umbilical points 
 for $2$-parameter families of surfaces \\
%{\red last ver}
}
\author[Y.~Kabata]{Yutaro Kabata}
\address[Y. ~Kabata]{Department of Mathematics,
Graduate School of Science,  Hokkaido University,
Sapporo 060-0810, Japan}
\email{kabata@mail.sci.hokudai.ac.jp}
\author[J. L. ~Deolindo Silva]{Jorge Luiz Deolindo Silva}
\address[J. L.~Deolindo-Silva]{Departamento de Ci\^encias Exatas e Educa\c{c}\~ao,
Universidade Federal de Santa Catarina-Blumenau-SC, Brazil}
\email{jorge.deolindo@ufsc.br}
\author[T.~Ohmoto]{Toru Ohmoto}
\address[T.~Ohmoto]{Department of Mathematics,
Faculty of Science,  Hokkaido University,
Sapporo 060-0810, Japan}
\email{ohmoto@math.sci.hokudai.ac.jp}
\begin{abstract} 
We determine local topological types of 
binary differential equations of asymptotic curves 
at parabolic and flat umbilical points 
for generic $2$-parameter families of surfaces in $\Proj^3$ 
by comparing our projective classification of Monge forms 
and classification of general BDE obtained by Tari and Oliver. 
In particular, generic bifurcations of the parabolic curve are classified. 
The flecnodal curve is also examined by direct computations, 
and we present new bifurcation diagrams in typical examples. 
\end{abstract}

%\begin{keyword}
%Projective differential geometry of surfaces, 
%binary differential equations, asymptotic curves,  parabolic curve, flecnodal curve, 
%singularities of smooth maps.
%\end{keyword}

%\end{frontmatter}
\maketitle

\setlength{\baselineskip}{16pt}

\section{Introduction}
{\it Binary differential equations} (BDE) widely appear  
in several geometric problems. 
A BDE in two variables $x, y$
has the form
\begin{equation}\label{bde}
 a(x,y)\, dy^2+2 b(x,y) \, dx dy + c(x,y) \, dx^2=0
 \end{equation}
with smooth functions $a, b, c$ of $x, y$. 
It is regarded as a smooth map $\R^2 \to \R^3$ 
assigning $(x, y) \mapsto (a, b, c)$ and
consider the $C^\infty$-topology on the mapping space. 
Put $\delta(x,y):=b(x,y)^2-a(x,y)c(x,y)$. 
If $\delta>0$,
the BDE locally defines two foliations which are transverse to each other. 
The {\it discriminant curve} is given by $\delta=0$, 
at which the integral curve of BDE generically has a cusp.  
Two germs of BDEs $F$ and $G$ are equivalent if
there is a local diffeomorphism in the $xy$-plane
sending the integral curves of $F$ to those of $G$. 
Also the topological equivalence is defined. 
There have been known several classification results for germs of (families of) BDE  
\cite{BT1, BT2, BFT, CBR, DR, Davydov1, Davydov2, Tari1, Tari2}.  
As a specific geometric setting, 
consider a surface locally given by $z=f(x,y)$; 
{\it asymptotic curves} are integral curves of the BDE 
\begin{equation}\label{abde}
f_{yy} \, dy^2 + 2 f_{xy} \, dx dy + f_{xx} \, dx^2=0
\end{equation} 
(called an {\em asymptotic BDE},  for short). 
Asymptotic BDEs form a thin subset of the space of general BDEs. 
The discriminant curve coincides with the {\it parabolic curve} in the surface theory; 
denote it by $\mathcal{P}$. 
Note that the above equivalence relation of BDE preserves the discriminant, 
but loses any information about inflection of integral curves, 
thus the theory of general BDE is really useful for analyzing parabolic curves 
but not for flecnodal curves at all. 

In this paper, 
we are interested in bifurcation phenomena of asymptotic BDE. 
In \cite{SKSO}, 
we studied generic $1$ and $2$-parameter families of surfaces 
in real projective $3$-space $\Proj^3$ and 
presented a classification of Monge forms under projective transformations 
in accordance with equisingularity types of central projection;   
indeed it is a natural extension of a well-known classification of jets of a {\it generic} surface 
given by Arnold-Platonova (cf. \cite{Arnold, Platonova, Landis}) 
and is related to a work of Uribe-Vargas on 
$1$-parameter bifurcations of parabolic and flecnodal curves \cite{UV}. 
There are in total 20 normal forms of parabolic and flat Monge forms 
up to codimension $4$ (Table \ref{main_table1} in \S 2). 
For each normal form, 
we will carefully check the criteria in topological classification of general BDE 
due to Davydov, Bruce and Tari \cite{BT1, BT2, BFT, Davydov1, Davydov2, Tari1, Tari2}, 
that determines 
the diffeomorphic/topological type of our asymptotic BDEs  
(Propositions \ref{BDE0}, \ref{BDE1} and \ref{BDE2}). 
It then turns out that 
asymptotic BDEs at parabolic points arising in generic $2$-parameter families of surfaces 
realize any generic types of IDE (implicit differential equation) of codimension $2$ 
classified in Tari \cite{Tari1}. 
The BDE at flat umbilical points is more remarkable; 
While our degenerate flat umbilic class $\Pi^f_2$ generically appears in $2$-parameter family, 
its asymptotic BDE is not equivalent to 
any type of BDE of codimension $2$ with $(a(0), b(0),c(0))=(0,0,0)$ classified in Tari \cite{Tari2}, 
but it is equivalent to the normal form 
\begin{equation}\label{oliver}
xdy^2+2ydxdy+x^2dx^2=0,
\end{equation}
which is actually one of types 
of codimension $3$ in the space of (general) BDEs 
 obtained by Oliver \cite{Oliver}  
(Remark \ref{rem_BDE2}). 

\begin{figure}
\centering
  \includegraphics[height=5cm]{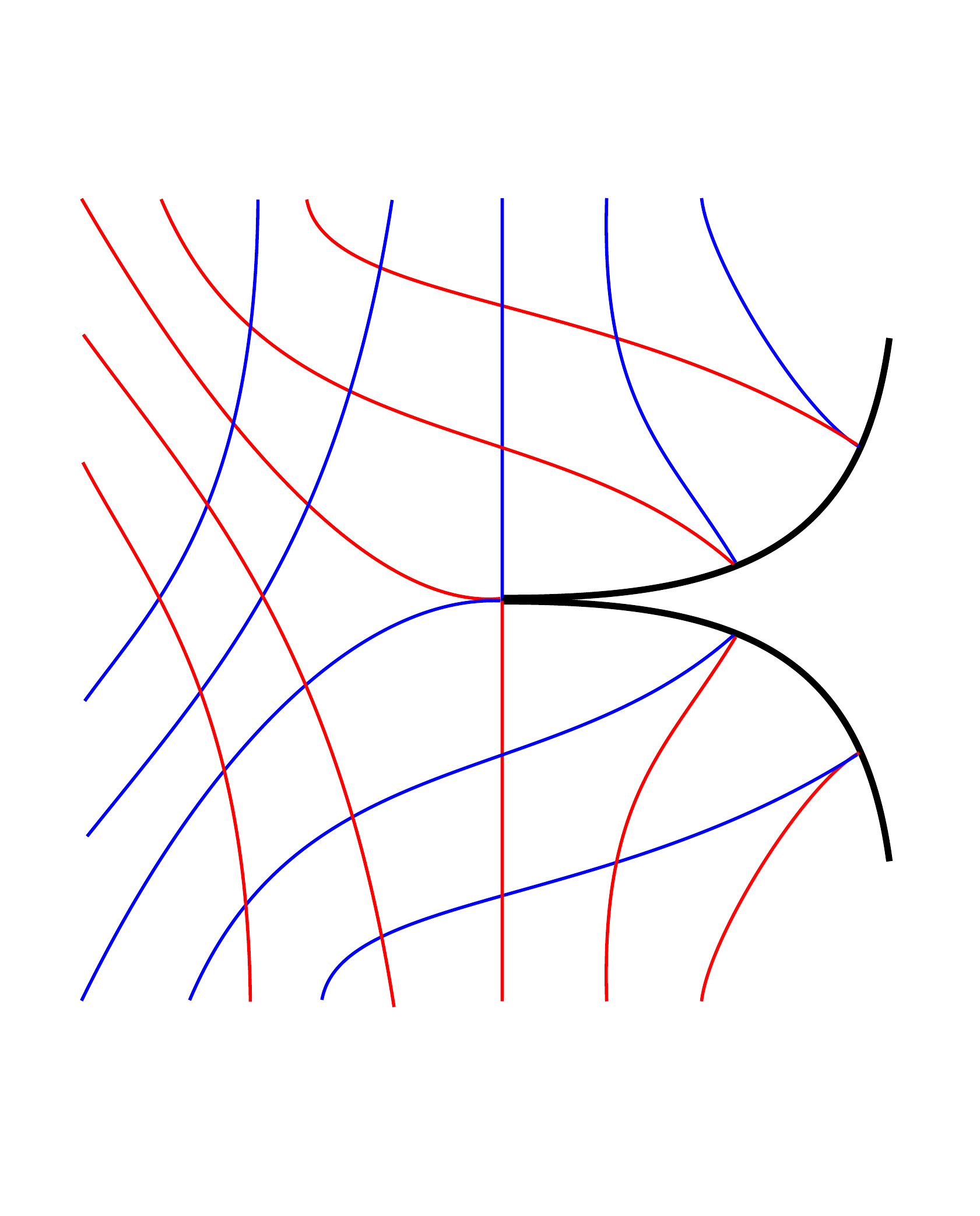}\;\;\; 
    \includegraphics[height=5cm]{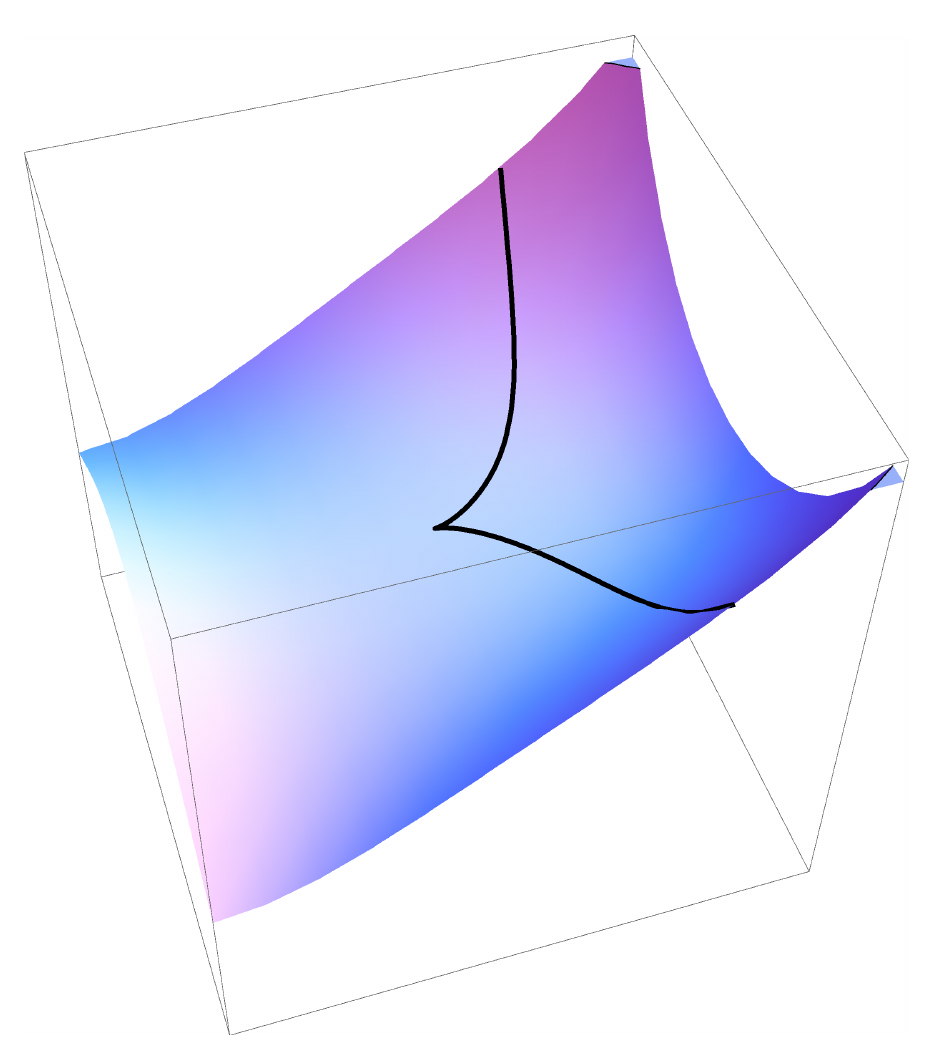}
\caption{BDE of type (\ref{oliver}) in \cite{Oliver} (left) 
and degenerate flat point ($D_5$) of type $\Pi_2^f$ in \cite{SKSO} (right). }
\label{oliver1}
\end{figure}

Next,  we find out the bifurcation diagrams for generic $2$-parameter families of surfaces. 
We show that obtained families of asymptotic BDEs 
are topologically versal (as families of general BDEs) 
in the sense of Tari \cite{Tari1}, except for the class $\Pi_2^f$ mentioned above. 
Then, bifurcations of the parabolic curve are read off from 
the bifurcation diagrams depicted in \cite{Tari1}. 
However, this is not useful for analyzing the flecnodal curve; 
for instance, unlike general BDE, 
the $A_3$-transition of asymptotic BDE at a point of type $\Pi_{v,1}^p$ creates 
a `figure-eight' flecnodal curve, as it was firstly { observed by  F. Aicardi  (Trieste, 1997) and 
A. Ortiz-Rodr\'iguez (Paris, 1999) through computer experiments} (cf. \cite{UV}). 
Therefore,  by direct computations, 
we examine the bifurcation of the flecnodal curve explicitly in examples. 
Also we present the bifurcation diagram of type  $\Pi^f_2$, 
which would completely be new in literature.

\subsection*{Acknowledgement}
The first and second authors thank organizers of the 14th Workshop on 
Real and Complex Singularities for giving them a nice opportunity to work together. 
The first author is supported by JSPS grant no.16J02200. 
The third author is partly supported by JSPS grants no.24340007 and 15K13452.

\section{Monge forms at parabolic and flat points} 
\begin{figure}
\centering
  \includegraphics[width=10cm]{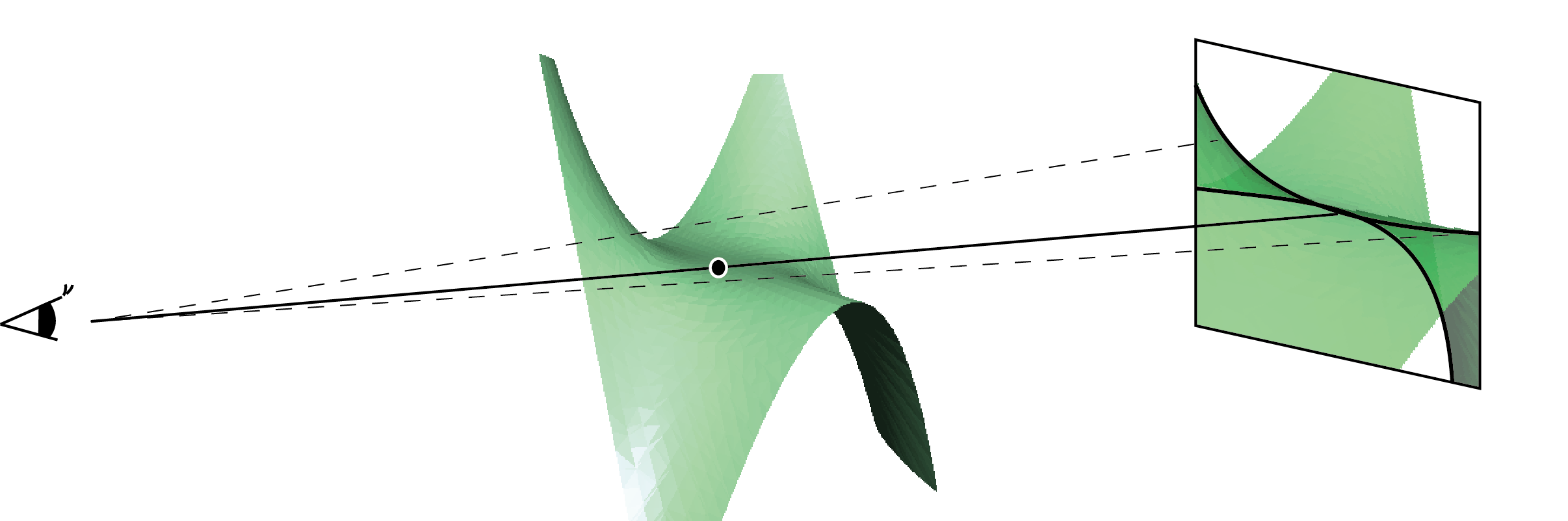}
  \caption{Central projection: 
  Parabolic points are characterized as points at which 
  the projection has the lips/beaks singularity or worse.}
  \label{beaks}
\end{figure}

We briefly recall our classification of Monge forms via projective transformations 
obtained in \cite{SKSO}. 
Take an affine chart $\R^3=\{[x:y:z:1]\} \subset \Proj^3$,
and consider germs of surfaces in $\R^3$ at the origin
given by Monge forms $z=f(x, y)$ with $f(0)=0$ and $df(0)=0$.
We say that two germs or jets of surfaces
are {\it projectively equivalent}
if there is a projective transformation on $\Proj^3$ sending one to the other.
Projective transformations preserving the origin and the $xy$-plane
form a $10$-dimensional subgroup of $PGL(4)$,
and it acts on the space $J=\mathfrak{m}_{x,y}^2/\mathfrak{m}_{x,y}^{k+1}$ 
of $k$-jets of Monge forms in a canonical way. 
In \cite{Platonova81, Platonova} (cf. \cite{Arnold, Landis}), 
Platonova studied a projectively invariant stratification of $J$ 
with codimension $\le 2$, 
and it has recently been extended by Kabata \cite{Kabata} systematically 
up to codimension $4$ 
so that each stratum is characterized 
by singularity types 
of central projections which the surface-germ possesses. 
Here, the {\it central projection} of a surface $M$ from a viewpoint $q \in \Proj^3$ 
means the restriction to $M$ of a canonical projection $\pi_q:\Proj^3-\{q\} \to \Proj^2$; 
at each point $p \in M$, 
the projection is locally described as a map-germ $\R^2, 0 \to \R^2, 0$ 
in local coordinates centered at $p \in M$ and $\pi_q(p) \in \Proj^2$, respectively, 
and consider its singularity type (the class up to {\it $\A$-equivalence}, i.e.,  
the equivalence relation of map-germs via 
natural actions of diffeomorphism-germs of the source and the target). 
The singularity type measures how the line contacts with $M$: 
from a point on non-asymptotic line, the projection is of type fold $\II_2: (y, x^2)$ 
($2$-point contact), 
and from a point of an asymptotic line, it is of cusp $\II_3:(y, x^3+xy)$ in general 
($3$-point contact), and 
a plenty of degenerate types of map-germs appear, 
which are not determined only by the contact order, e.g.  
the parabolic curve $\mathcal{P}$ is formed by points at which 
the projection has the beaks/lips singularity $\I_2: (y, x^3\pm x^2y)$ or worse (Figure \ref{beaks}). 

The normal forms of parabolic and flat umbilical Monge forms are listed in 
Table \ref{main_table1} below (\cite{SKSO, Platonova}), 
where $k$ is the order of jets, 
$\mbox{cod}$ is the codimension of strata, 
and the last column $\mbox{proj.}$ means singularity types of 
central projection of the surface at the origin 
from viewpoints on the asymptotic line 
(the type in bracket indicates a more degenerate singularity type 
of projection from some isolated viewpoint specially chosen on the line). 

Let $M \subset \Proj^3$ be a non-singular surface. 
Suppose that an open subset $U \subset M\cap \R^3$ is parametrized by 
a graph $z=f(x,y)$ of a function. 
Since we are working in projective geometry, %the situation is much simpler; 
we may define 
the Monge form at $p$ by the $k$-jet of $f$ at $p$ off the linear term $j^1f(p)$. 
Then the {\it Monge-Taylor map} $\theta: U \to J$  is locally defined 
with respect to this affine chart.  
Take an open cover of $M$ by such affine charts so that 
$M$ is locally given by graphs of functions. 
By a standard transversality argument, we easily see that 
any class  $\Pi^*_{*,*}$ with codimension $\le 2$ appears 
for a generic embedding $M \to \Proj^3$ (\cite{Platonova81}), and 
any class with $\mbox{cod}=k>2$ 
appears in a generic family 
of embeddings $M \times V \to \Proj^3$ ($V \subset \R^{k-2}$).

\begin{table}[h]
{%\small 
$$
\begin{array}{l | l | c c | l  }
\mbox{class}  & \mbox{normal form} & k & \mbox{cod}
& \mbox{proj. }  \\
\hline
\hline
\Pi^p_{\I,1} & y^2+x^3 + xy^3+\alpha x^4 & 4 & 1 &
\I_2 \, (\I_3)
\\
\Pi^p_{\I,2} & y^2+x^3\pm xy^4+\alpha x^4+\beta y^5 + x^2\phi_3 & 5 & 2 &  \I_2 \, (\I_4)
\\
\Pi^p_{c,1} & y^2+x^2 y + \alpha x^4 \;\;\; (\alpha\not=0, \frac{1}{4}) & 4 & 2 &  \III_2\, (\III_3)
\\
\hline
\Pi^p_{c,2}& y^2+x^2 y + \frac{1}{4} x^4 + \alpha x^5+ y \phi_4 \; (\alpha\not=0) & 5 & 3 &  \III_2 
\\
\Pi^p_{c,4} & y^2+ x^2 y + x^5 + y\phi_4 & 5 & 3 & \IV_1
\\
\Pi^p_{\I,3} & y^2+x^3+ xy^5+\alpha x^4+\phi  & 6 & 3 &  \I_2\,  (\I_{5})\\
\Pi^p_{v,1}& y^2\pm x^4 +\alpha x^3y+\beta x^2y^2 \;\; (\beta\not=\pm\frac{3}{8}\alpha^2) & 4 & 3 & \V_1\,  (\VI)
\\
\Pi^f_{1} & xy^2\pm x^3 + \alpha x^3 y + \beta y^4 & 4  & 3 & \I_2^\pm,  \I_3  (\I_4)  \\
\hline 
\Pi^p_{c,3} & y^2+x^2 y + \frac{1}{4} x^4 + y \phi_4  & 5 & 4 &  \III_3\, (\III_{4})
\\
\Pi^p_{c,5} & y^2+ x^2 y \pm x^6+y(\phi_4+\phi_5)    & 6 & 4 &  \IV_2
\\
\Pi^p_{\I,4} &y^2+x^3+\alpha x^4+\phi & 6 & 4 &  \I_2\,  (\I_{6})\\
\Pi^p_{v,2}& y^2 \pm x^4+ \alpha x^3y \pm \frac{3}{8}\alpha^2 x^2y^2  & 4 & 4 & \V_1\,  (\VI_1)\\
\Pi^p_{v,3} & y^2+ x^5 + y(\phi_3+\phi_4)  & 5 & 4 & \V_2\, (\VI_2)
\\
\Pi^f_{2} & xy^2 + x^4 \pm y^4+\alpha x^3 y  & 4 & 4 & \I_2^- (\III) \\
\end{array}
$$
}
\caption{\small  Monge forms at parabolic and flat points are obtained in 
\cite{Arnold, Platonova, Landis} for $\mbox{cod}=1,2$ and \cite{SKSO} for $\mbox{cod}=3,4$.  
In the list,  $\alpha, \beta, \cdots$ are leading moduli parameters, 
$\phi_r$ denotes generic homogeneous polynomials of degree $r$ 
and $\phi=\beta y^5 + \gamma y^6 + x^2(\phi_3 + \phi_4)$. 
Double-sign $\pm$ corresponds in same order for each of 
$\Pi^p_{v,1}$ or $\Pi^p_{v,2}$. 
}
\label{main_table1}
\end{table}

\section{Binary differential equations}

\subsection{General BDE}
\begin{figure}
\centering
  \includegraphics[width=10cm]{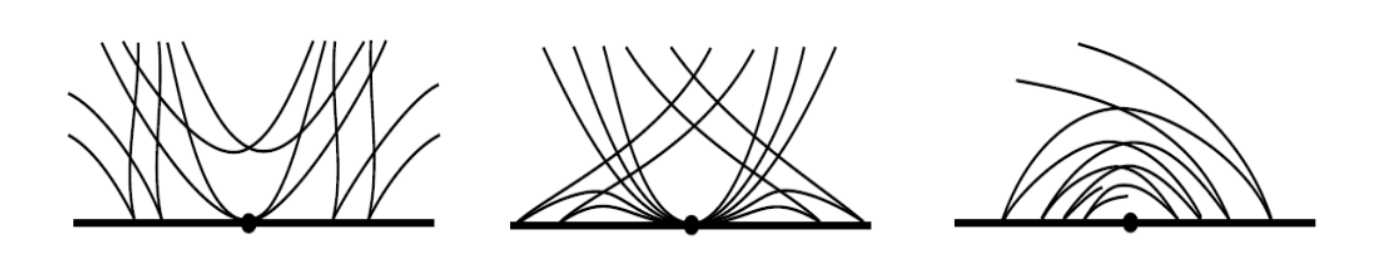}\\
  \caption{Folded singularities: saddle (left), node (center) and focus (right)}\label{fig1}
\end{figure}

One can separate BDE (\ref{bde}) into two cases.
The first case occurs
when the functions $a, b, c$ do not vanish at the origin at once,
then the BDE is just an implicit differential equation (IDE).
The second case is that all the coefficients of BDE vanish at the origin.
Stable topological models of BDEs belong to the first case;
 it arises when the discriminant is smooth (or empty).
If the unique direction at any point of the discriminant is transverse to it
(i.e. integral curves form a family of cusps),
then the BDE is stable and smoothly equivalent to
$dy^2+xdx^2=0$, 
that was classically known in Cibrario \cite{CBR} and also Dara \cite{DR}.
If the unique direction is tangent to the discriminant,
then the BDE is stable and smoothly equivalent to
$dy^2+(-y+\lambda x^2)dx^2=0$
with $\lambda\neq 0,\frac{1}{16}$, that was shown in Davydov \cite{Davydov1, Davydov2}; 
the corresponding point in the plane is called a \emph{folded singularity} 
-- more precisely, 
a \emph{folded saddle} if $\lambda<0$, 
a \emph{folded node} if $0<\lambda<\frac{1}{16}$ 
and a \emph{folded focus} if $\frac{1}{16}<\lambda$, see  Figure \ref{fig1}.

In both cases,
the topological classification of
generic $1$ and $2$-parameter families of BDEs have been established
in Bruce-Fletcher-Tari \cite{BT1, BFT} and Tari \cite{Tari1, Tari2}, respectively.
We will use those results later.
Besides, we need
a generic $3$-parameter family  of BDE studied in Oliver \cite{Oliver}.

\subsection{BDE of asymptotic curves} 
We are concerned with (degenerate) parabolic points 
and flat umbilic points of a surface. 
In fact, asymptotic BDE is 
intimately related to the singularity type of the Monge form. 
The parabolic curve can be seen as the locus
where the Monge form has $A_{\geq2}$-singularities;
when the Monge form has a $A_3^\pm$-singularity,
the surface has a cusp of Gauss,
which corresponds to the class $\Pi^p_{c,1}$ --
in this case the asymptotic BDE has a folded saddle singularity 
(resp. a folded node or focus singularity) 
if the Monge form has a singularity of type $A_3^-$ (resp. $A_3^+$). 
The transitions in $1$-parameter families occur generically
in three ways at the following singularities of the Monge form:
non-versal $A_3$,  $A_4$ and  $D_4$ (flat umbilic) \cite{BFT}.
For $2$-parameter families,
$A_3,$  $A_4$,  $A_5$ and $D_5$ singularities of the Monge form 
generically appear.
Below, the Monge form is written by 
$$f(x,y)=\sum_{2\le i+j} c_{ij}\, x^i y^j, $$ 
and the $k$-jet $j^kf(0)$ is assumed to coincide with the normal form as 
 in Table \ref{main_table1} for each class. 

\begin{prop} \label{BDE0}
The following classes in Table \ref{main_table1}
correspond to structurally stable types of BDE given  in \cite{CBR, DR, Davydov1}. 
\begin{enumerate}
\item[$(\Pi^p_{\I,k})$]
$(1 \le k\le 4)$.  
The parabolic curve is smooth
and the unique direction defined by $\delta=0$ is transverse to the curve; 
the asymptotic BDE is smoothly equivalent to 
$$dy^2+x\, dx^2=0.$$
\item[$(\Pi^p_{c,k})$]
$(k=1,4,5)$. 
The parabolic curve is smooth
and the unique direction defined by  $\delta=0$ is tangent to the curve; 
the asymptotic BDE is smoothly equivalent to
$$dy^2+(-y+ \lambda x^2)dx^2=0$$
with $\lambda=6(c_{40}-\frac{1}{4})\not=0$,
where $c_{40}$ is the coefficient of $x^4$ in the normal form.
\end{enumerate}
\end{prop}

\proof The results follow from the comments in \S 3.1 above.
In  second case, i.e., $j^4f= y^2+x^2y+c_{40}x^4$, 
the $2$-jet of the asymptotic BDE is transformed to the above form 
via $x=\bar{x}$ and $y=-\frac{1}{2}\bar{x}^2-\bar{y}$.
\qed

\begin{rem}\label{rem_BDE0}\upshape 
As $c_{40}=0$ in the normal forms of classes $\Pi^p_{c,4}$ and $\Pi^p_{c,5}$, 
we see $\lambda=-\frac{3}{2}<0$,
thus the asymptotic BDE has a folded saddle at the origin. 
The {\it folded saddle-node bifurcation} (cf. Fig.2 in \cite{Tari1}) occurs at $\lambda=0$. 
That is the case of $c_{40}=\frac{1}{4}$,
that corresponds to the classes $\Pi^p_{c,k}\; (k=2,3)$ dealt below. 
Notice that another exceptional value $\lambda=\frac{1}{16}$ does not relate to
our classification of Monge forms given by projection-types (Table \ref{main_table1}). 
That is, the {\it  folded node-focus bifurcation} of asymptotic BDE occurs
within the same class $\Pi^p_{c,1}$ (cf. Fig.3 in \cite{Tari1}) and 
$\lambda=\frac{1}{16}$ makes a condition on coefficients of order greater than $4$ 
of the normal form, that is independent from 
the geometry of central projection of the surface. 
We should remark that 
$\Pi^p_{c,4}$ and $\Pi^p_{c,5}$ cause 
$1$ and $2$-parameter bifurcations of the flecnodal curve, respectively. 
For instance, during a $1$-parameter bifurcation of type $\Pi^p_{c,4}$, 
a butterfly point moves on the flecnodal curve 
and passes through this degenerate cusp of Gauss 
at the bifurcation moment, see \cite[\S 4]{SKSO} and \cite{UV}. 
\end{rem}

\begin{rem}\label{rem_BDE1}\upshape 
At elliptic points in an smooth surface in $\mathbb R^3$ there is a unique pair of conjugate directions for which the
included angle (i.e the angle between these directions) is minimal. These directions are called {\it characteristic directions} and are determined in terms of the coefficients of the first and second fundamental forms. Theses directions are not preserved via projective transformations, but at a cusp of Gauss, Oliver in \cite{Oliver2} shows that the characteristic directions are invariant under projective transformations. We can use the normal form  $\Pi^p_{c,1}$ to obtain the BDE associated to characteristic directions; it is indeed smoothly equivalent to $dy^2+(-y+\lambda x^2)dx^2=0$ with $\lambda=-6c_{40}+\frac{3}{2}\not=0$. The configurations of asymptotic and characteristic curves at a cusp of Gauss are given in  \cite{Oliver2}. 
\end{rem}

\begin{prop}\label{BDE1}
The following classes 
correspond to some topological types of BDE with codimension $1$.
\begin{enumerate}
\item[$(\Pi^p_{v,1})$] 
The Monge form has
an $A_3$-singularity at the origin, at which
the parabolic curve has a Morse singularity;
the asymptotic BDE is topologically equivalent to 
the non-versal $A_3^\pm$-transitions with
Morse type 1  in \cite{BFT}
$$dy^2+(\pm x^2 \pm y^2) dx^2=0.$$
\item[$(\Pi^p_{c,2})$]
The Monge form has
an $A_4$-singularity at the origin, at which
the parabolic curve is smooth;
the asymptotic BDE is topologically equivalent to 
the well-folded saddle-node type in \cite{BFT, Davydov2}
$$dy^2+(-y + x^3) dx^2=0,$$
provided the coefficient of $x^5$ in the normal form $c_{50}\neq 0$.
\item[$(\Pi^f_{1})$\,\;]
The Monge form has
a $D_4^\pm$-singularity at the origin, at which
the parabolic curve has a Morse singularity;
the asymptotic BDE is topologically equivalent to 
the bifurcation of star/$1$-saddle types in \cite{BT1}
\begin{eqnarray*}
D_4^+: && ydy^2-2x dxdy -y dx^2=0 \;\; \mbox{\rm (star)};\\
D_4^-: &&ydy^2+2x dxdy +y dx^2=0 \;\; \mbox{\rm ($1$-saddle)}.
\end{eqnarray*}
\end{enumerate}
\end{prop}

\proof 
In \cite{SKSO, Kabata}, the class $\Pi^p_{v,1}$ is explicitly described as follows. 
Let $z_0=y^2+c_{20}x^2+\sum_{i+j\ge 3} c_{ij}x^iy^j \in J$. 
Then, $z_0$ is projectively equivalent to $\Pi^p_{v,1}$ if and only if 
$$c_{20}=c_{30}=c_{21}=0, \; \;  c_{40}\neq0, \;\;
S:=3 c_{31}^2 + 8 c_{40}(c_{12}^2  - c_{22})\not=0.$$ 
In fact, exactly the same condition appears in \cite[p.501, Case 1]{BFT} 
as the condition for $A_3^\pm$-transition: 
$S\not=0$ means the $2$-jet $j^2\delta(0)$ is non-degenrate (ibid),  
thus the normal form follows from Theorem 2.7 (and Prop. 4.1)  in \cite{BFT}.
For the class $\Pi^p_{c,2}$,  
$z_0$ (of the above form) is projectively equivalent to $\Pi^p_{c,2}$ 
if and only if 
$$c_{20}=c_{30}=B=0, \; c_{40}\not=0, \; A\not=0,$$
with $B:=c_{21}^2c_{40}-4c_{40}^2$ and 
$A:=c_{21}^2c_{50}+4c_{12}c_{40}^2-2c_{21}c_{31}c_{40}$. 
This condition is the same as the one for $A_4$-transition in  \cite[p.502, Case 2]{BFT}, 
and then the normal form of BDE is obtained  
(\cite[Prop. 4.2]{BFT}, also see \cite{Davydov2}). 
For the class $\Pi^f_{1}$,
the asymptotic BDE is given in \cite[Cor. 5.3]{BT1}:
indeed, for our normal form of $\Pi^f_{1}$,
the parabolic curve is defined by
$3 x^2 -  y^2 + 18 \beta x y^2+ \cdots =0$,
hence it has a node at the origin for arbitrary $c_{31}=\alpha, c_{04}=\beta$.
\qed

\begin{prop}\label{BDE2}
The following classes 
correspond to some topological types of BDE with codimension $\ge 2$. 
\begin{enumerate}
\item[$(\Pi^p_{v,2})$] 
The Monge form has
an $A_3$-singularity at the origin,
at which the parabolic curve has a cusp singularity;
the asymptotic BDE is topologically equivalent to the cusp type  in \cite{Tari1}
$$dy^2+(\pm x^2+y^3) dx^2=0,  $$
provided 
$C_1
:=\mp 5c_{50}c_{31}^3+12c_{41}c_{31}^2\mp 24c_{32}c_{31}+32c_{23}\neq  0
$.
\item[$(\Pi^p_{v,3})$] 
The Monge form has
an $A_4$-singularity at the origin,
at which the parabolic curve has a Morse singularity;
the asymptotic BDE is topologically equivalent to 
the non-transversal Morse type  in \cite{Tari1}
$$dy^2+(xy+x^3) dx^2=0$$
provided $C_2:=c_{31}\neq 0$.
\item[$(\Pi^p_{c,3})$] 
The Monge form has an $A_5$-singularity at the origin,
at which the parabolic curve is smooth;
the asymptotic BDE is topologically equivalent to 
the folded degenerate elementary type  in \cite{Tari1}
$$dy^2+(-y \pm x^4) dx^2=0,$$
provided $C_3:=c_{60}-\frac{1}{2}c_{41}\neq 0$.
\item[($\Pi^f_{2})$\,\;] 
The Monge form has
a $D_5$-singularity at the origin, at which
the parabolic curve has a cusp singularity;
the asymptotic BDE is topologically equivalent to 
a cusp type 2 in \cite{Oliver}
$$xdy^2+2ydxdy+x^2dx^2=0.$$
\end{enumerate}
\end{prop}

\proof 
For each of the first three classes,
the claim follows from Proposition 4.1 and Theorem 1.1 of Tari \cite{Tari1}. 
Let $S, A, B$ be as in the proof of Proposition \ref{BDE1}. 
As shown in \cite{SKSO}, 
the condition of $z_0=y^2+c_{20}x^2+o(2)$ to be equivalent to $\Pi_{v,2}^p$ is given by 
$$c_{20}=c_{30}=c_{21}=S=0,\;\; c_{40}\neq0$$
which is entirely the same as the condition of (iii) in \cite[p.156]{Tari1} 
($C_1$ is given by $C$ in the bottom of that page). 
Also the condition for  $\Pi_{v,3}^p$ is given by 
$$c_{20}=c_{30}=c_{21}=c_{40}=0,$$
and that for $\Pi_{c,3}^p$ is given by 
$$c_{20}=c_{30}=B=A=0.$$
The same conditions can be found in (ii) and  (i) in \cite[p.156]{Tari1} respectively 
($C_2, C_3$ are given by $c_3$  in (ii) and $A$  in (i) in {\it ibid.}). 
Hence, those asymptotic BDEs are equivalent to 
the normal forms presented in Theorem 1.1 in \cite{Tari1}. 
For the last class  $\Pi^f_{2}$,
the $1$-jet of the asymptotic BDE is given by $j^1(a,b,c)(0)=(2x,2y,0)$ and
the parabolic curve is defined by $- 4 y^2+24 x^3 +\cdots =0$,
namely it has a cusp at the origin. 
Thus the corresponding BDE is equivalent to 
one of types described in Theorem 3.4 of  \cite{Oliver}.
\qed

\begin{rem}\label{rem_BDE2}
\upshape
In Tari \cite{Tari2} and Oliver \cite{Oliver},
BDE with the discriminant having a cusp are classified. 
Notice that the BDE for $\Pi^f_{2}$ in Proposition \ref{BDE2} is equivalent to 
one of `type 2' in Oliver's classification \cite{Oliver} 
which appears in a generic $3$-parameter family of general BDE's, 
while the type appears in a generic $2$-parameter family of asymptotic BDE 
by our classification of Monge forms \cite{SKSO}. 
This is not surprising, for asymptotic BDEs form a thin subset of the space of all BDEs. 
In fact, it is shown in \cite[Prop.2.1(2)]{Tari2} that 
for general BDE with a cusp, the $1$-jet is reduced 
by linear changes of coordinates and multiplication by non-zero constants to 
$j^1(a,b,c)(0)=(x, \pm  y+ \alpha x, 0)$ ($\alpha \in \R$),  
while in our case it is reduced to the particular form $j^1(a,b,c)(0)=(x, y, 0)$ 
as seen in the proof of Proposition \ref{BDE2}. 
This infers the gap of codimensions caused by two different classifications. 
It would be interesting to find a deeper understanding of the geometry 
of asymptotic BDE at a flat umbilic point. 
\end{rem}

\begin{rem}
\upshape
In \cite{SKSO, Kabata}, our classification of Monge forms has been achieved 
in accordance with singularity types of central projections, or say almost equivalently, 
the contact of the surface with lines, while 
Tari \cite{Tari1} described types of asymptotic BDEs in terms of 
singularities of height functions and singularities of parabolic curves, 
that reflects the contact of the surface with planes. 
These two different approaches lead to the same conditions $C_1, C_2, C_3\not=0$. 
That should be explained by using 
 a duality between the contact of lines and the contact of planes. 
\end{rem}

\section{Families of Monge forms and BDEs}

In Propositions \ref{BDE1} and \ref{BDE2},   
we have compared our Monge forms in \cite{SKSO} 
and types in classification of general BDE given by Tari \cite{Tari1}.   
In this section, we compare families of Monge forms and families of BDE.  
Given an $s$-parameter family 
$f(x,y,\lambda)\; (=f_{\lambda}(x,y)): U \times \R^s \to \R$  $(U \subset \R^2)$, 
we define a family of Monge-Taylor maps 
$$\theta: U \times \R^{s} \to J, \quad 
\theta(p,\bu):=j^kf_{\lambda}(p).$$ 
Below, 
for each class in Table \ref{main_table1}, 
we take a family of Monge forms whose Monge-Talyor map $\theta$ is transverse to 
the corresponding stratum in the jet space $J$ (Table \ref{table2}). 
We show that the associated family of asymptotic BDEs 
is topologically versal in the sense of Tari \cite{Tari1}. 

\subsection{Transverse families of Monge forms}
For instance, 
recall the cases of $\Pi_{v,k}^p$ as in 
Propositions \ref{BDE1} and \ref{BDE2}. 
Write $z=\sum_{2\le i+j\le k} c_{ij}x^iy^j \in J$ as before, 
and regard $c_{ij}$ as coordinates of $J$.  
The locus of parabolic Monge forms in $J$ is defined by $c_{20}c_{02}-c_{11}^2=0$, 
thus the tangent space to the locus at $z_0=y^2+o(2)$ is defined by 
the $1$-form $dc_{20}=0$ 
in $T_{z_0}J=J$. 

Let $z_0=y^2+o(2)$ be of type $\Pi_{v,1}^p$. 
By the local defining equation of the stratum, 
its tangent space at $z_0$  is given 
by linear equations 
$$dc_{20}=dc_{30}=dc_{21}=0\quad \mbox{on} \;\; T_{z_0}J.$$
In particular, take $z_0=y^2\pm x^4 +\alpha x^3y+\beta x^2y^2+o(4)$ 
with $S(z_0)=3\alpha^2-8\beta\not=0$ 
and $f: U \to \R$ a representative of it;  $z_0=j^4f(0)$.  
The Monge-Taylor map $\theta: U \to J$ sends  
$p \mapsto j^4f(p)$ off the linear term. 
Then, the image $d\theta(T_0U)$ and $\rd/\rd c_{20}$ 
span the normal space (i.e the quotient of $T_{z_0}J$ via the tangent space of the stratum) 
and thus the $1$-parameter family $f(x,y,t)=f(x,y)+tx^2$ 
induces a family of Monge-Taylor maps $\theta: U\times \R \to J$ 
being transverse to the stratum at the origin $(0,0)$. 

For the class $\Pi_{v,2}^p$, 
let 
$$\textstyle 
f(x,y)=y^2\pm x^4+\alpha x^3y\pm \frac{3}{8}\alpha^2 x^2y^2+\sum_{i+j=5}c_{ij}x^iy^j+o(5).$$ 
Then the tangent space of the stratum at $z_0=j^4f(0)$ is defined by 
$$dc_{20}=dc_{30}=dc_{21}=dS=0\quad \mbox{on} \;\; T_{z_0}J.$$
Since 
$dS=6c_{31}dc_{31}+8(c_{12}^2-c_{22})dc_{40}+16c_{12}c_{40}dc_{12}-8c_{40}dc_{22}$, 
we have $dS=6\alpha dc_{31}\mp 3 \alpha^2 dc_{40} -8dc_{22}$ at $z_0$. 
Then the condition that
$$\textstyle 
\mbox{rank}\; [\; dc_{20}, dc_{30}, dc_{21}, dS\;]^T
\left[d\theta(\frac{\rd}{\rd x})\; d\theta(\frac{\rd}{\rd y})\right]=2$$ 
is written down as  
$$C_1:=\mp 5c_{50}c_{31}^3+12c_{41}c_{31}^2
\mp 24c_{32}c_{31}+32c_{23}\neq  0,$$
that is exactly the same condition required in Proposition \ref{BDE2} (\cite[p.156]{Tari1}). 
Then we can easily find a desired $2$-parameter deformation of $f$; 
for instance, when $c_{23}\not=0$ and other $c_{50}=\cdots=c_{05}=0$ 
(then $C_1\not=0$), 
we may take $f(x,u,t,u)=f(x,y)+t x^2+u x^2y$. 

Also for other cases in Table \ref{main_table1}, 
any representative $f(x,y)$ of the normal form ($k$-jet of Monge form) 
admits a deformation whose Monge-Taylor map is transverse to the stratum, 
provided Taylor coefficients of $f$ of higher order ($> k$) are chosen 
to be appropriately generic, if necessary.  

In Table 2, 
we collect examples of such families of Monge forms 
deforming the normal forms in Table \ref{main_table1}. 
Here we omit the stable case dealt in Proposition \ref{BDE0}.

 \begin{table} 
$$
\begin{array}{l | l l}
\mbox{class} & \mbox{family}\\
\hline\hline
\Pi^p_{c,2}& y^2+x^2 y + \frac{1}{4} x^4 + \alpha x^5%+ y \phi_4 
+ t x^3  &(\alpha\not=0)
\\
\Pi^p_{c,3}(\pm) & y^2+x^2 y + \frac{1}{4} x^4 +\gamma x^4y %+ y \phi_4  +\phi_6
 + t x^3 + u x^4 & (\gamma \lessgtr 0)
\\
%\Pi^p_{c,4} & y^2+ x^2y + x^5 + y\phi_4 + t x^3 &\\
%\Pi^p_{c,5} & y^2+ x^2 y \pm x^6+y(\phi_4+\phi_5)  + t x^3+u x^4 &\\
%\hline
%\Pi^p_{\I,3} & y^2+x^3+xy^5+\alpha x^4+\phi  + t x y^3 &\\
%\Pi^p_{\I,4} &y^2+x^3+ xy^6+\alpha x^4+\phi +t xy^3+u x y^4 &\\
\hline
\Pi^p_{v,1}(\pm, +) & y^2+ x^4 +\alpha x^3y+\beta x^2y^2 + t x^2  & (\beta> \frac{3}{8}\alpha^2)
\\
\Pi^p_{v,1}(\pm, -) & y^2-x^4 +\alpha x^3y+\beta x^2y^2 + t x^2  & (\beta< \frac{3}{8}\alpha^2)
\\
\Pi^p_{v,2}(\pm) & y^2 \pm x^4 + \alpha x^3y \pm \frac{3}{8}\alpha^2 x^2y^2 
+ \gamma x^2y^3+ t x^2+u x^2y & (\gamma\not=0)
\\
\Pi^p_{v,3} & y^2+ x^5 + \gamma x^3y + t x^2+u x^2y & (\gamma\not=0)
\\
\hline
\Pi^f_{1}(\pm) & xy^2\pm x^3 + \alpha x^3 y + \beta y^4 + t x^2
\\
\Pi^f_{2}(\pm) & xy^2 + x^4 \pm y^4+\alpha x^3 y + t x^2+ux^3 
\end{array}
$$
\caption{\small  Examples of families of Monge forms (with parameters $t, u$) 
whose Monge-Taylor maps are transverse to the strata. 
}
\label{table2}
\end{table}

\subsection{Versal families of BDE}
As a general theory, 
the germ of  BDE (\ref{bde}) with $a(0)\not=0$ 
is canonically transformed to the germ of an IDE 
$$p^2+\frac{a(x,y)c(x,y)-b(x,y)^2}{a(x,y)^2}=0 \qquad \left(p=\frac{dy}{dx}\right)$$ 
by a simple coordinate change $\bar{x}=x$ and $\bar{y}=y+\int_0^x \frac{b(s,y)}{a(s,y)}ds$. 
For an IDE, $p^2+\varphi(x,y)=0$, 
moreover for a family of IDEs, $p^2+\varphi(x,y,t,u)=0$, 
there are known useful criteria for detecting its genericity, by which 
generic classifications of IDE and that of families of IDEs have been achieved, 
see Tari \cite{Tari1} for the detail (also \cite{BFT}). 
Given a deformation of a parabolic Monge form $y^2+o(2)$, 
we obtain a family of IDEs of asymptotic curves, and apply Tari's criteria to it. 

Now we check criteria of asymptotic IDE and BDE 
for the following four examples of families of Monge forms in Table \ref{table2}: 

\begin{itemize}
\item[($\Pi_{v,2}^p$)]
$p^2\pm 6 x^2 + \gamma y^3+ u y +t - \frac{3}{2}u^2 x^2 - 3 \gamma t x^2 y=0$ 
$(\gamma\not=0)$ (we set $\alpha=0$ for simplicity). 
It is just an IDE of cusp type in the sense of \cite{Tari1}. 
Check a  criterion in Proposition 3.5 (ii) of \cite{Tari1}: 
$$
\left|
\begin{array}{cc}
\varphi_t(0) & \varphi_{ty}(0) \\
\varphi_u (0) & \varphi_{uy}(0) 
\end{array}\right|
=-1\not=0, 
$$
thus, by Theorem 3.6 of \cite{Tari1}, 
our family of asymptotic BDE is fiber topologically equivalent to
$$dy^2+(\pm x^2+y^3+uy +t) dx^2=0. $$
\item[($\Pi_{v,3}^p$)]
$\textstyle p^2+3 \gamma x y+10 x^3+t  + u y - \frac{3}{2} u^2 x^2  - 5 \gamma u x^3=0$ 
$(\gamma\not=0)$. 
It is of non-transverse Morse singularity type. 
Check a  criterion in Proposition 3.3 (ii) of \cite{Tari1}: 
$$
 \left|
\begin{array}{ccc}
0&1&6\\
\varphi_t(0) & \varphi_{ty}(0) &  \varphi_{txx}(0)\\
\varphi_u(0) & \varphi_{uy}(0) &  \varphi_{uxx}(0)\\
\end{array}\right|
=-6\not=0, 
$$
thus, by Theorem 3.4 of \cite{Tari1}, 
our family is fiber topologically equivalent to
$$dy^2+(xy+x^3+ux^2+t) dx^2=0.$$
\item[($\Pi_{c,3}^p$)] 
$\textstyle p^2
+ y \mp \frac{15}{2}\gamma x^4+3 t x + 6 u x^2  \pm 6 \gamma x^2 y=0$ 
$(\gamma\not=0)$. 
It is of folded degenerate elementary singularity type. 
Check a  criterion in Proposition 3.1 (ii) of \cite{Tari1}: 
$$
 \left|
\begin{array}{cc}
\varphi_{tx}(0) & \varphi_{txx}(0) \\
\varphi_{ux} (0) & \varphi_{uxx}(0) 
\end{array}\right|
=-36\not=0, 
$$
thus, by Theorem 3.2 of \cite{Tari1}, 
our family  is fiber topologically equivalent to
$$dy^2+(-y\pm x^4+ux^2+tx) dx^2=0.$$
\end{itemize}

\begin{itemize}
\item[($\Pi_{f,1}^p$)] 
This is not the case of IDE and indeed it is a $1$-parameter family of BDE, 
thus we refer to Example 4.1 in \cite{BT2}.  
Consider 
$F=(12\beta y^2+2x)p^2+2(3\alpha x^2+2y)p+(6\alpha xy+2t+6x)=0$. 
In this case the linear part of $F$ provides 
all the topological information about the family of BDEs (\cite{BT1,BT2}). 
Check a versality criterion in Proposition 2.1 of \cite{BT2}: 
$$
 \left|
\begin{array}{ccc}
2 & 0&0 \\
0&2  & 0\\
6&0&2 
\end{array}\right|
=8\not=0, 
$$
thus we can reduce the 1-jet of the BDE  to the form $(y+t)dy^2\pm 2xdxdy \pm y dx^2=0$. 
Combining Theorem 3.5  and Example 4.1 in \cite{BT2} 
with $\phi(p)=(F_x+pF_y)(0,0,0,p)=6p^2+6$, 
our family  is fiber topologically equivalent to
\begin{eqnarray*}
(y+t)dy^2-2 x dxdy - y dx^2&=&0,\\
(y+t)dy^2+2x dxdy + y dx^2&=&0.
\end{eqnarray*}
\end{itemize}

Also for families of Monge forms of type  $\Pi_{v,1}^p$, $\Pi_{c,2}^p$, 
in Table \ref{table2}, it can be seen that the families of asymptotic BDE 
are respectively equivalent to 
\begin{eqnarray*}
dy^2+(\pm x^2+\pm y^2+t) dx^2&=&0,  \mbox{(see \cite{BFT})}\\
dy^2+(-y+x^3+tx) dx^2&=&0\;  \mbox{(see \cite{Tari1})}.
\end{eqnarray*}

Bifurcation diagrams of generic $2$-parameter families of IDEs 
have clearly been depicted in Tari \cite{Tari1, Tari2}. 
Therefore, we can deduce from those figures 
the bifurcation diagrams of parabolic curves 
for generic $2$-parameter families of parabolic Monge forms. 
In the next section, we compute the bifurcation of flecnodal curve 
at parabolic and flat umbilical points.

\section{Bifurcation diagrams for $2$-parameter families of surfaces}

\subsection{Flecnodal curve}
A point of a surface in $\Proj^3$ is {\it flecnodal} 
if an asymptotic line at that point has more than $3$-point contact with the surface. 
The closure of the set of such points is called the {\it flecnodal curve}, 
denoted by $\mathcal{S}$, which is an important characteristic of the surface; 
$\mathcal{S}$ lies on the hyperbolic domain and 
meets the parabolic curve $\mathcal{P}$ 
at (ordinary or degenerate) cusps of Gauss $\Pi_{c,*}^p$. 
A flecnodal point is characterized as the point 
at which the projection along an asymptotic line is of type 
the swallowtail singularity $\II_4: (y, x^4+yx)$ or worse. 
From this fact, 
a local defining equation of $\mathcal{S}$ is obtained in a very neat way  
\cite{Saji, Kabata}.  
Suppose that the origin $0\in \R^3$ is 
a flecnodal point of a surface $z=f(x,y)$ such that 
the $x$-axis is the asymptotic line at $0$. 
For deforming the line, one has $2$-dimensional freedom, thus 
the projection along the $x$-axis, $(x,y)\mapsto (y, f(x,y))$, has 
 a $2$-parameter deformation 
$$F_{v,w}(x,y) = (y-vx, f(x,y)-wx).$$ 
Let $\lambda=0$ be the equation 
defining the singular point set (contour generator) of $F_{v,w}$ 
and $\eta$ be a vector field on a neighborhood of the origin in $\R^2$ 
which spans $\ker  dF_{v,w}$ where $\lambda=0$, 
i.e. 
$$\textstyle 
\lambda(x,y,v,w):=\det dF_{v,w}(x,y) \;\; \mbox{and} \;\;  
\eta(x,y,v,w):=\frac{\rd}{\rd x}+v\frac{\rd}{\rd y}.$$ 
Then the swallowtail singularity of $F_{v,w}$, and thus the curve $\mathcal{S}$, 
is characterized by three equations 
$$\lambda=\eta\lambda=\eta\eta\lambda=0.$$ 
By $\lambda=0$, $w$ is always solved.  
Eliminating $v$ by the last two equations, we obtain  
an equation of variables $x, y$ or parametrizations $x=x(v)$, $y=y(v)$, 
which defines $\mathcal{S}$ around $(x,y)=(0,0)$. 
Furthermore, there generically appear some isolated points in the curve $\mathcal{S}$ 
at which an asymptotic line has $4$-point contact with the surface, 
i.e. the projection $F_{v,w}$ admits the butterfly singularity  $\II_5: (y, x^5+yx)$, 
so we call such a point a {\it butterfly point} for short. 
It is defined by an additional equation 
$\eta\eta\eta\lambda=0$ on $\mathcal{S}$.

The parabolic curve $\mathcal{P}$ is obtained as 
the locus of points where the singular point set of 
the projection $F_{v,w}$ is not smooth, i.e., 
$$\textstyle
\lambda=\frac{\rd \lambda}{\rd x}=\frac{\rd \lambda}{\rd y}=0.$$ 
In Figures below, 
$\mathcal{P}$ is drawn in black and $\mathcal{S}$ is in gray.

\subsection{$1$-parameter bifurcations} 
For three classes of codimension $1$ in Proposition \ref{BDE1},  
we confirm bifurcations of curves $\mathcal{P}$ and $\mathcal{S}$ as depicted in \cite{UV}  
by direct computations using families of Monge forms in Table \ref{table2}. 
In \cite{UV} 
bifurcations of $\mathcal{S}$ at hyperbolic points are also classified 
and positive/negative flecnodal points are considered. 
We do not enter the full scope of the classification of hyperbolic points, 
since we focus mainly on bifurcations at parabolic and flat umbilical points. 

\

\begin{figure}[p]
\centering
    \includegraphics[width=11.5cm]{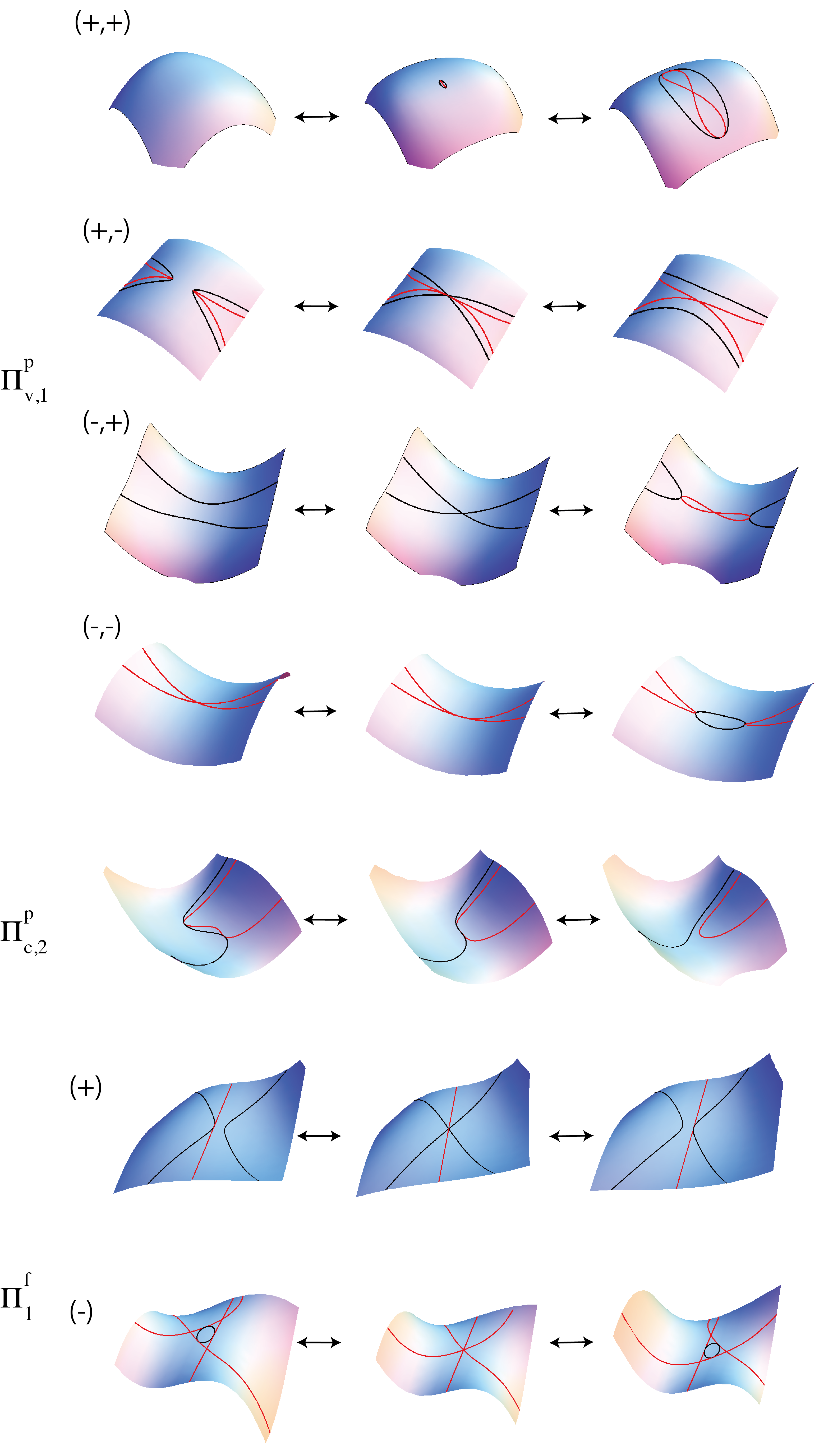}
  \caption{$1$-parameter bifurcations of $\Pi_{v,1}^p$, $\Pi_{c,2}^p$ and $\Pi_1^f$ 
  (\cite{UV}, \cite{SKSO}).}
  \label{1para}
\end{figure}

\t
$\bullet \;\; (\Pi^p_{v,1}(\pm,\pm))$ 
As seen in Proposition \ref{BDE1}, at a point of type $\Pi^p_{v,1}$, 
non-versal $A_3^\pm$-transition of BDE occurs (cf. \cite{SKSO}). 
According to sign difference of coefficients in the normal form, there are four types. 
Among them, there are two types so that 
the flecnodal curve is created/canceled when passing through the point. 
Unexpectedly, the flecnodal curve has the form of `figure-eight',  { as mentioned in Introduction. }
Obviously, the Proposition \ref{BDE1} 
does not help anything for understanding the appearance of the figure eight curve, 
because the equivalence of BDE does not preserve inflections of integral curves. 
Let us confirm this fact by 
a direct computation using the normal form of $\Pi^p_{v,1}$. 
Let 
$$f(x,y, t)=y^2+ x^4 + x^2y^2 + t x^2.$$
Solving equations $\lambda=\eta\lambda=\eta\eta\lambda=0$, 
we have 
$$\textstyle
(x,y)=
\left(\mp \frac{v(2+v^2)\sqrt{-t-v^2}}{{\sqrt{2}(-2+v^2)\sqrt{2+v^2-v^4}}}, 
\pm \frac{(2+v^2)\sqrt{-t-v^2}}{\sqrt{2(2+v^2-v^4)}}\right)$$
with $t\le -v^2$ and $|v| \ll 1$, which parametrizes   
part of the flecnodal curve $\mathcal{S}$ 
sitting in the half planes, $y\ge0$ and $y\le 0$. 
The parabolic curve $\mathcal{P}$ is given by 
$t + (6+t) x^2 + y^2  + 6 x^4 - 3 x^2 y^2=0$. 
As $t$ varies, an elliptic Morse bifurcation of the parabolic curve occurs 
and the created flecnodal curve has the form of figure-eight 
as depicted in $(+, +)$, Figure \ref{1para}. 
No butterfly point appears on $\mathcal{S}$, 
since $\eta\eta\eta\lambda\not=0$ for $(x,y)$ near the origin.  

Also for the form $f(x,y, t,u)=y^2- x^4 + x^2y^2 + t x^2$,  
it has a hyperbolic Morse bifurcation of the parabolic curve 
and the figure-eight flecnodal curve also arises as depicted in $(-, +)$. 
In the other two cases, $f(x,y, t,u)=y^2\pm x^4 - x^2y^2 + t x^2$, 
the curves bifurcate as depicted in  $(\pm, -)$.

\

\t
$\bullet \;\; (\Pi^p_{c,2})$ 
In this case, $A_4$-transition of asymptotic BDE occurs; during this process, 
a pair of cusps of Gauss (tangential points of $\mathcal{P}$ and $\mathcal{S}$) 
is created/canceled. Take 
$$\textstyle 
f(x,y, t)=y^2 + x^2 y + \frac{1}{4} x^4 + x^5  + t x^3$$
and project it along the $x$-axis. 
Then $\mathcal{P}$ is a cubic curve: 
$2y + 20 x^3+x^2 +6 tx=0$, and 
$\mathcal{S}$ is given by 
$y+100 x^4+10 x^3+ (\frac{1}{2} + 20 t) x^2 + 3 t x +t^2=0$.

\

\t
$\bullet \;\; (\Pi^f_1(\pm))$ 
Take 
$$\textstyle
f_\pm (x,y,t)=x y^2 \pm x^3 + x^3 y  + t x^2.$$
As a test, we consider the projection along the $x$-axis (and $y$-axis). 
In case of $f_+$, 
$\mathcal{P}$ is given by $12 x^2 - 9 x^4 - 4 y^2+4 t x =0$, 
and $\mathcal{S}$ is given by $x=0$. 
In case of $f_{-}$, 
$\mathcal{P}$ is given by 
$12 x^2 + 9 x^4 + 4 y^2-4 t x =0$, and $\mathcal{S}$ has three branches 
consisting of 
the $y$-axis and two smooth curves having an intersection point moving along the $x$-axis.

\subsection{$2$-parameter bifurcations} 
For each of four classes of codimension $2$ in Proposition \ref{BDE2}, 
we draw a new picture of the bifurcation diagram using the family in Table \ref{table2}. 
The bifurcation of $\mathcal{P}$ can be indeed read off from 
Figure 8, Figure 6 and Figure 4  in Tari \cite{Tari1}, respectively.  
We compute the curve $\mathcal{S}$ 
and add some new branches to Tari's figures. 

\

\begin{figure}
\centering
  \includegraphics[width=11.5cm]{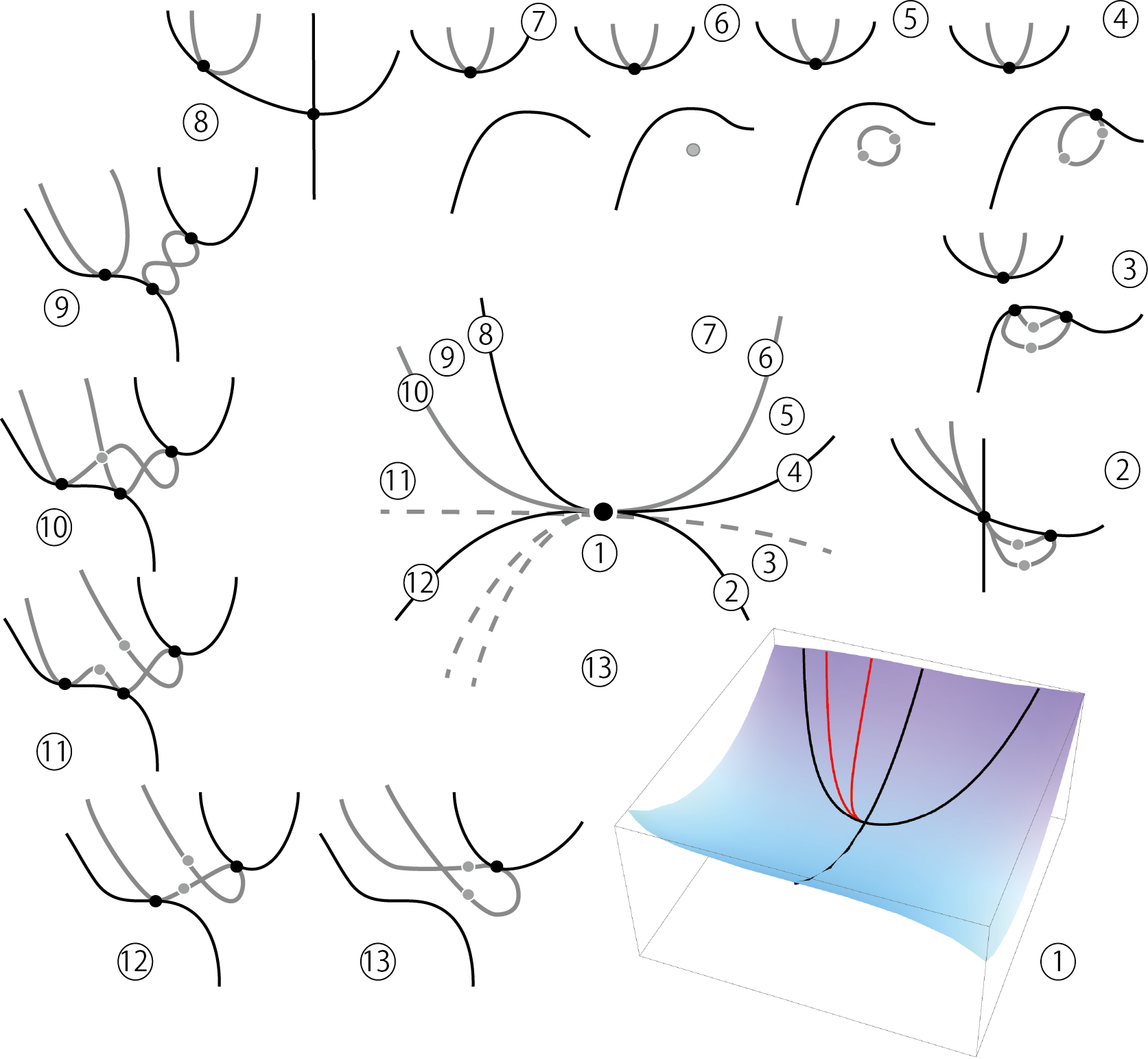}
  \caption{Bifurcation of $\Pi_{v,3}^p$. 
The diagram consists of branches named by 
  $\Pi_{v,1}^p$ (no.2, 8),   $\Pi_{c,2}^p$ (no.4, 12),  $\Pi_{3,5}^h$  (no.6, 10),  
  $\Pi_{4,5}^h$ (between no.12-13) and $\Pi_{c,4}^p$  (between no.3-4). 
  }
    \label{Piv3}
\end{figure}

\t
$\bullet \;\; (\Pi_{v,3}^p)$  
Consider the family of Monge forms  
$$f(x,y, t,u)=y^2+ x^5 + x^3y + t x^2+u x^2y,$$
and then we have Figure \ref{Piv3} below (cf. Figure 8 in \cite{Tari1}). 
The parabolic curve $\mathcal{P}$ is given by the equation 
$p(x,y,t,u)=12 x y+ 40 x^3- 9 x^4+4 t + 4 u y - 4 u^2 x^2  - 12 u x^3  =0$. 
When $(t,u)=(0,0)$, 
$\mathcal{P}$ consists of the $y$-axis and the smooth curve 
$y=-\frac{10}{3} x^2 + \frac{3}{4} x^3$; hence, $\mathcal{P}$ has a node at the origin. 
Solving $p=p_x=p_y=0$, 
we see that 
a hyperbolic Morse bifurcation of $\mathcal{P}$ 
occurs at a point of type $\Pi_{v, 1}^p$ 
when the parameter comes across a smooth curve 
$t=\frac{1}{27}u^3 (10+\frac{3}{4} u)$ (no.2 and 8 in Figure \ref{Piv3}). 
The $A_4$-transition at which two cusps of Gauss are created/canceled appears 
along a smooth curve on $ut$-space, that is the $\Pi_{c,2}^p$-locus  (no.4 and 12).  
On the other hand, 
solving equations $\lambda=\eta\lambda=\eta\eta\lambda=0$ for $y$, 
$\mathcal{S}$ is expressed by two branches 
$$\textstyle
y=-\frac{1}{2} (u^3 +7 u^2 x + 20 x^2 + 18 u x^2 + 18 x^3 
\pm (u + 3 x)h(x,y,t,u))$$
where $h=
\sqrt{-4 t + u^4 + 8 u^3 x + (40 u + 28 u^2) x^2 + (80 + 48 u) x^3 + 36 x^4}$. 
In particular, 
the defining equation is written as 
$y^2 +100 x^4 + 20 x^2 y + 18 x^3 y  + u^3 y + t (u + 3 x)^2 + 
 u^2 (-10 x^3 + 7 x y) - 6 u (5 x^4 - 3 x^2 y)=0$. 
If $(t,u)=(0,0)$, 
$\mathcal{S}$ has a $5/2$-cusp at the origin which is tangent to the $x$-axis. 
By $\eta\eta\eta \lambda=0$, 
one can find some points of $\mathcal{S}$ at which 
the butterfly singularity appears in the projection; 
such isolated butterfly points are traced in Figure \ref{Piv3}. 
Two dotted branches between no.12 and 13 
indicate the bifurcation of class  $\Pi_{4,5}^h:\;xy+x^4+y^5+\alpha xy^3+\beta x^3y+x\phi_4$ 
in \cite[\S 5]{SKSO}, 
where a butterfly point  
passes through a double point of $\mathcal{S}$,  
and a dotted curve between no.3 and 4 
corresponds to the class $\Pi_{c,4}^p$ as noted in  Remark \ref{rem_BDE0},  
where a butterfly point passes through a degenerate cusp of Gauss. 
A lengthy computation shows that 
the butterfly point degenerates  
into the class  $\Pi_{3,5}^h:\;xy+x^3+y^5+\alpha xy^3+x\phi_4$ with $\alpha=0$ in \cite[\S 4]{SKSO}, 
when the parameter $(t,u)$ lies on 
a smooth curve $t=\frac{1}{4}u^4+o(4)$; 
an elliptic Morse bifurcation of $\mathcal{S}$ appears on one half branch (no.6), 
and a hyperbolic Morse bifurcation appears on the other branch (no.10).

\

\begin{figure}
\centering
    \includegraphics[width=8cm]{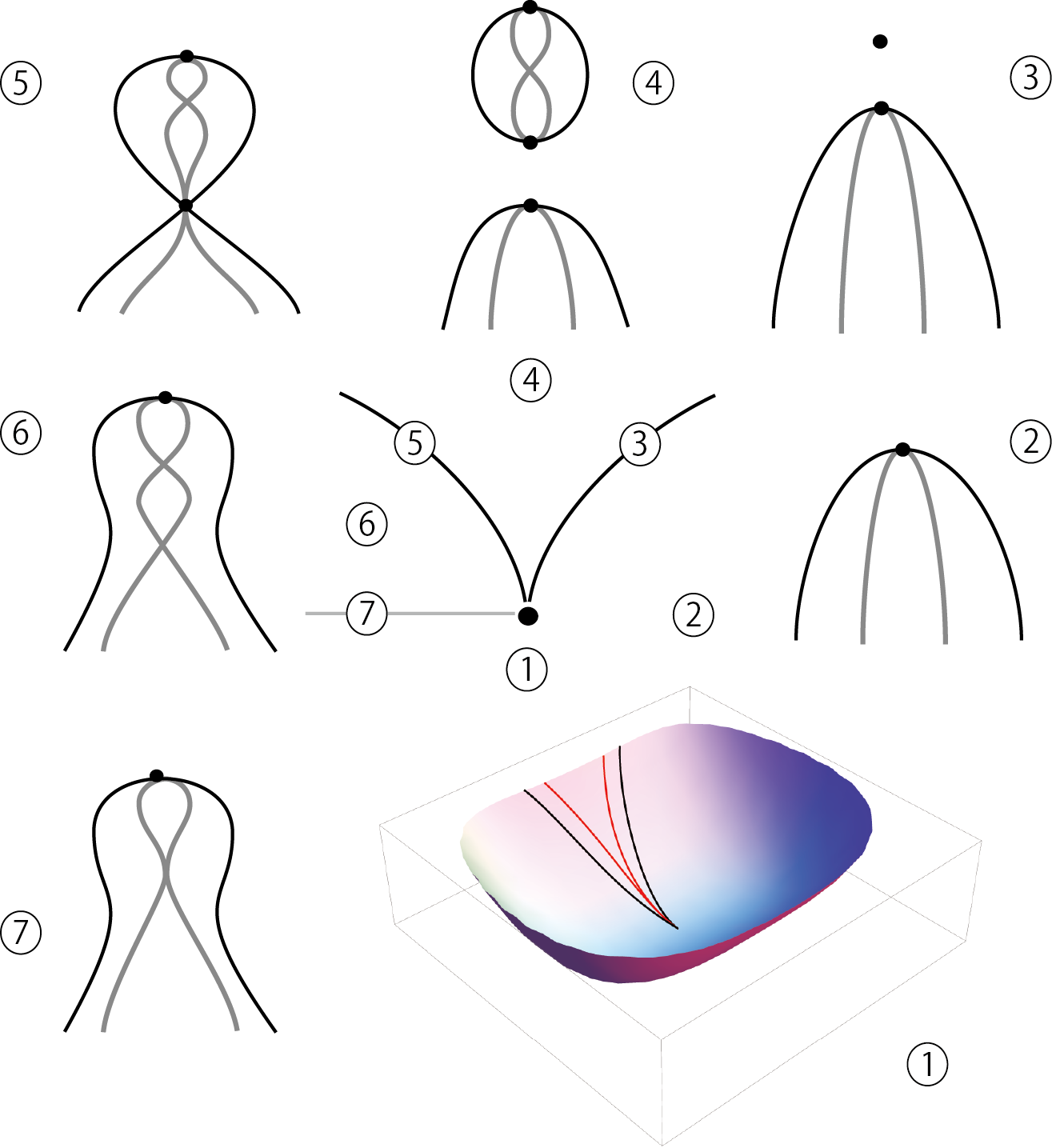}
    \caption{Bifurcations of $\Pi_{v,2}^p(+)$.}
    \label{Piv2+}
\end{figure}

\begin{figure}
\centering
  \includegraphics[width=8cm]{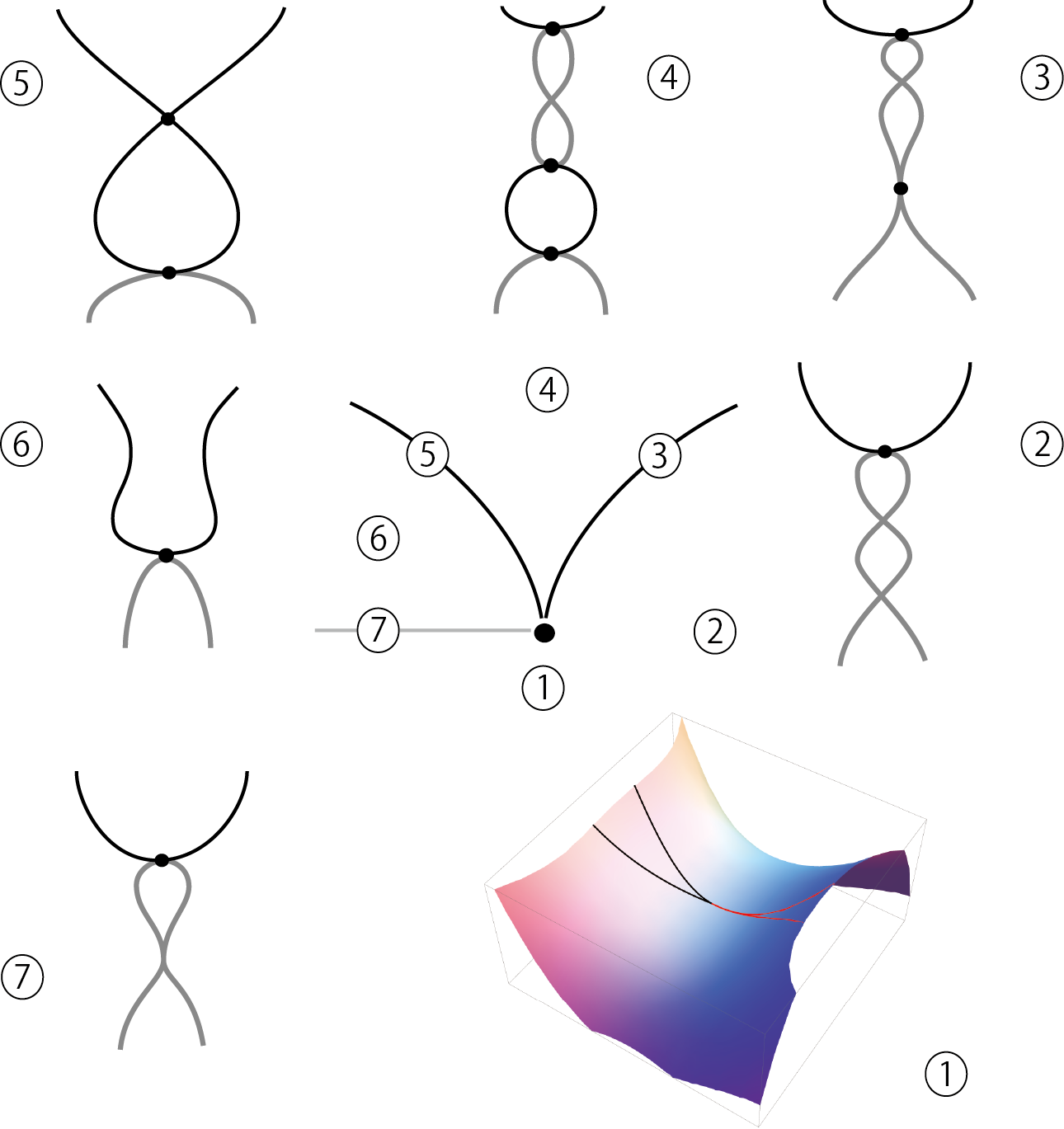}
    \caption{Bifurcations of $\Pi_{v,2}^p(-)$.}
    \label{Piv2-}
\end{figure}

\t
$\bullet \;\; (\Pi_{v,2}^p(\pm))$ 
Let 
$$\textstyle
f(x,y, t,u)=y^2 \pm x^4 + x^2 y^3  + t x^2 + u x^2 y.$$
In entirely the same way,  we have Figures 
 \ref{Piv2+} and \ref{Piv2-} (cf. Figure 6 in \cite{Tari1}). 
When $(t,u)=(0,0)$, 
$\mathcal{P}$ is defined by 
$x^2 \pm \frac{1}{6}y^3+ 3 x^4 y \mp x^2 y^4 =0$, and 
$\mathcal{S}$ is  $16 x^2 + 9 y^7\mp 66 x^2 y^4+o(7)=0$. 
No butterfly point occurs for $\eta\eta\eta\lambda\not=0$. 


\

\begin{figure}
\centering
  \includegraphics[width=9cm]{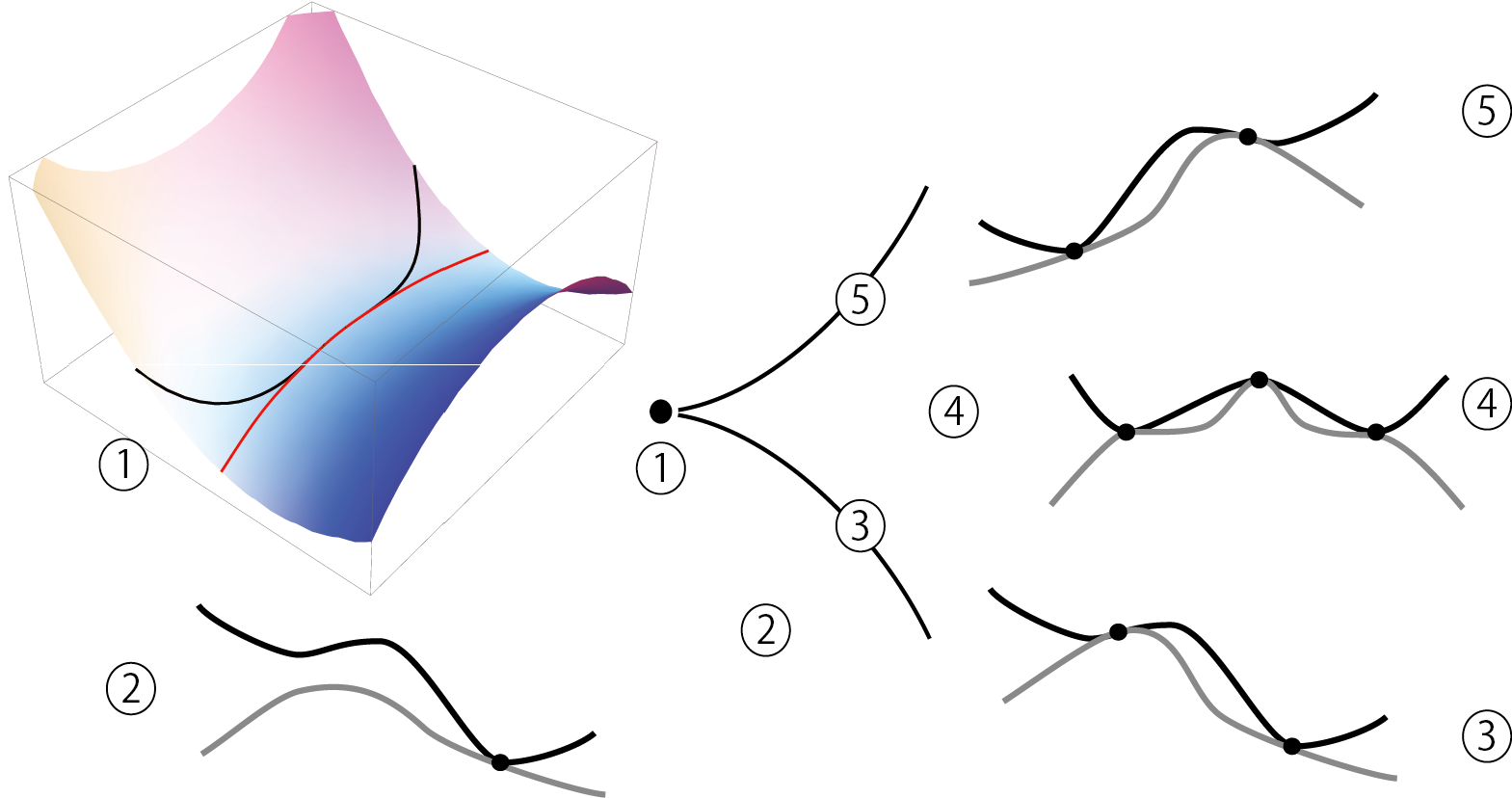}
  \caption{Bifurcations of $\Pi_{c,3}^p(+)$.}
  \label{Pic3+}
\end{figure}

\begin{figure}
\centering
  \includegraphics[width=9cm]{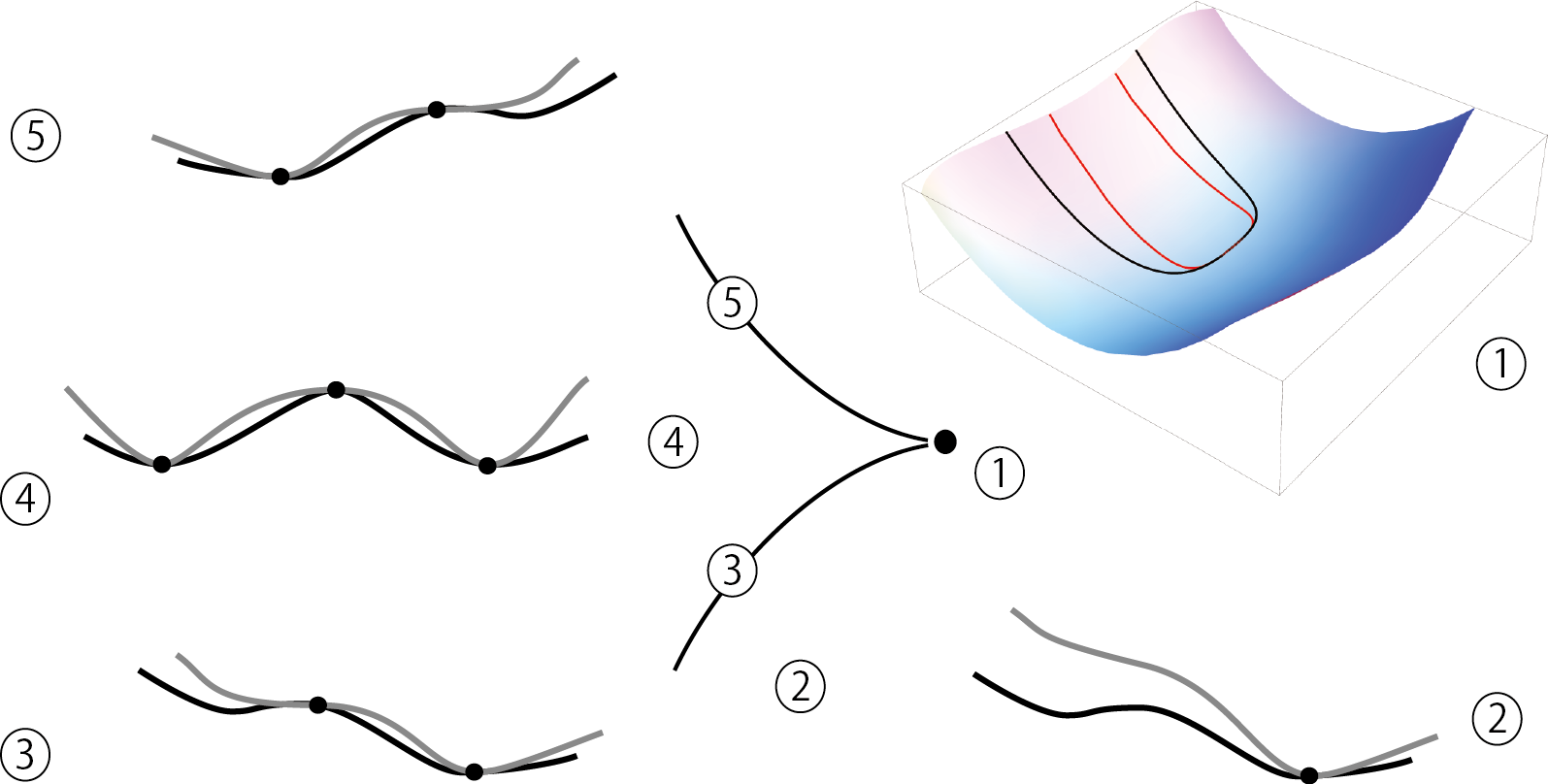}
  \caption{Bifurcations of $\Pi_{c,3}^p(-)$.}
  \label{Pic3-}
\end{figure}

\t
$\bullet \;\; (\Pi_{c,3}^p(\pm))$ 
Consider the family of Monge forms 
$$\textstyle
f(x,y, t,u)=y^2 + x^2 y + \frac{1}{4} x^4 \pm x^4 y  + t x^3 + u x^4.$$
In entirely the same way,  we have Figures \ref{Pic3+} and \ref{Pic3-} 
(cf. Figure 4  in \cite{Tari1}). 
When $(t,u)=(0,0)$, 
$\mathcal{P}$ and $\mathcal{S}$ are  
$y=-\frac{1}{2}x^2\pm7 x^4 - 38 x^6 +o(6)$ and 
$y=-\frac{1}{2}x^2 \pm 7 x^4 - 138 x^6+o(6)$, respectively. 
No butterfly point occurs for $\eta\eta\eta\lambda\not=0$. 


\

$\bullet \;\; (\Pi_2^f(\pm))$   
This case must be new, since it is not versal as a family of BDE. Take 
$$f(x,y, t,u)=xy^2+ x^4 \pm y^4 + t x^2+u x^3.$$ 
First, consider the case $(+)$ here (the coefficient of $y^4$ is $+1$). 
Any tangent line is asymptotic, so choose as a test the projection along the $y$-axis 
and its deformation: $(x,y)\mapsto (x-vy, f(x,y)-wy)$. 
The result is as follows; see Figure \ref{Pif2+} below. 
It is easy to find the equation of $\mathcal{P}$: 
$p(x,y,t,u)=6 x^3 - y^2 +t x + 3 u x^2 + 6 t y^2 + 18 u x y^2 + 36 x^2 y^2=0$. 
Hence, when $(t,u)=(0,0)$, $\mathcal{P}$ has an ordinary cusp at the origin. 
Solve $p=p_x=p_y=0$, then 
we have an equation, $t(32t-12u^2)=0$, 
that defines the locus of $(t,u)$ where $\mathcal{P}$ has a singularity at some $(x,y)$.  
When one comes across the component $t=0$ (no.6 and 9 in Figure \ref{Pif2+}) 
the bifurcation of type $\Pi_1^f$ occurs, 
and around the curve $32t=12u^2$  (no.4 and 11) 
the bifurcation of type  $\Pi_{v,1}^p$ appears.  
Solving $\lambda=\eta\lambda=\eta\eta\lambda=0$, 
$x,y$ can be parametrized by $v$, 
and then we can draw the curve $\mathcal{S}$ 
(we remark that $v$ may go to $\infty$, that means that 
the asymptotic line at such a point of $\mathcal{S}$ is close to the $x$-axis; 
then it is better for analyzing $\mathcal{S}$ to switch to the projection along the $x$-axis). 
In case of $(t,u)=(0,0)$, $\mathcal{S}$ has an ordinary cusp 
together with two smooth components passing through the origin. 
The bifurcation of $\mathcal{S}$ occurs as follows: 
At no.3, a tacnode bifurcation (self-tangency of two branches) appears 
(that class is denoted by $\Pi_{4,4}^h:\;xy+x^4\pm y^4+\alpha xy^3+\beta x^3y$ in \cite{SKSO}). 
Also between no.7 and 8 (also between no.13 and 2), 
two tacnode bifurcations arise successively, 
and from no.12 to no.13, 
there appear two events of type $\Pi_{c,2}^p$,
at each of which an $A_4$-transition appears, i.e. 
two cusps of Gauss 
(at which $\mathcal{P}$ and $\mathcal{S}$ are tangent) are canceled. 
There must be two branches between no.7-8 (also no.13-2, no.12-13) 
when taking general coefficients of order $5$ in the normal form $f$ 
(for the above particular form,  it is observed 
by the symmetry of $y\leftrightarrow -y$ that 
these two branches duplicate). 
No butterfly point occurs for $\eta\eta\eta\lambda\not=0$ near the origin. 
The case $\Pi_2^f(-)$ is slightly simpler than $\Pi_2^f(+)$ (Figure \ref{Pif2-}).

\begin{figure}
\centering
\includegraphics[width=11cm]{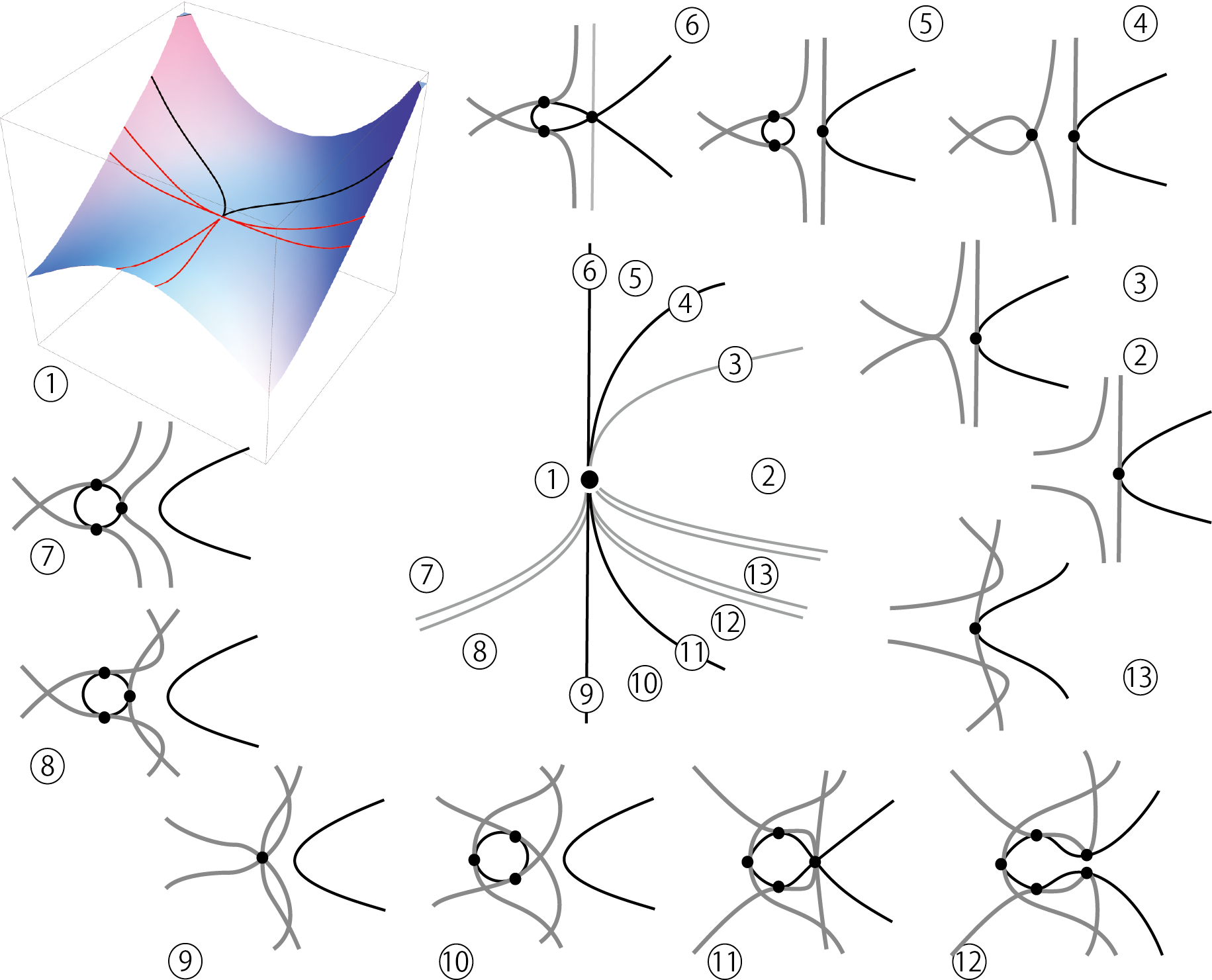}
\caption{Bifurcation of $\Pi_2^f(+)$.}
\label{Pif2+}
\end{figure}

\vspace{1cm}

\begin{figure}
\centering
\includegraphics[width=10cm]{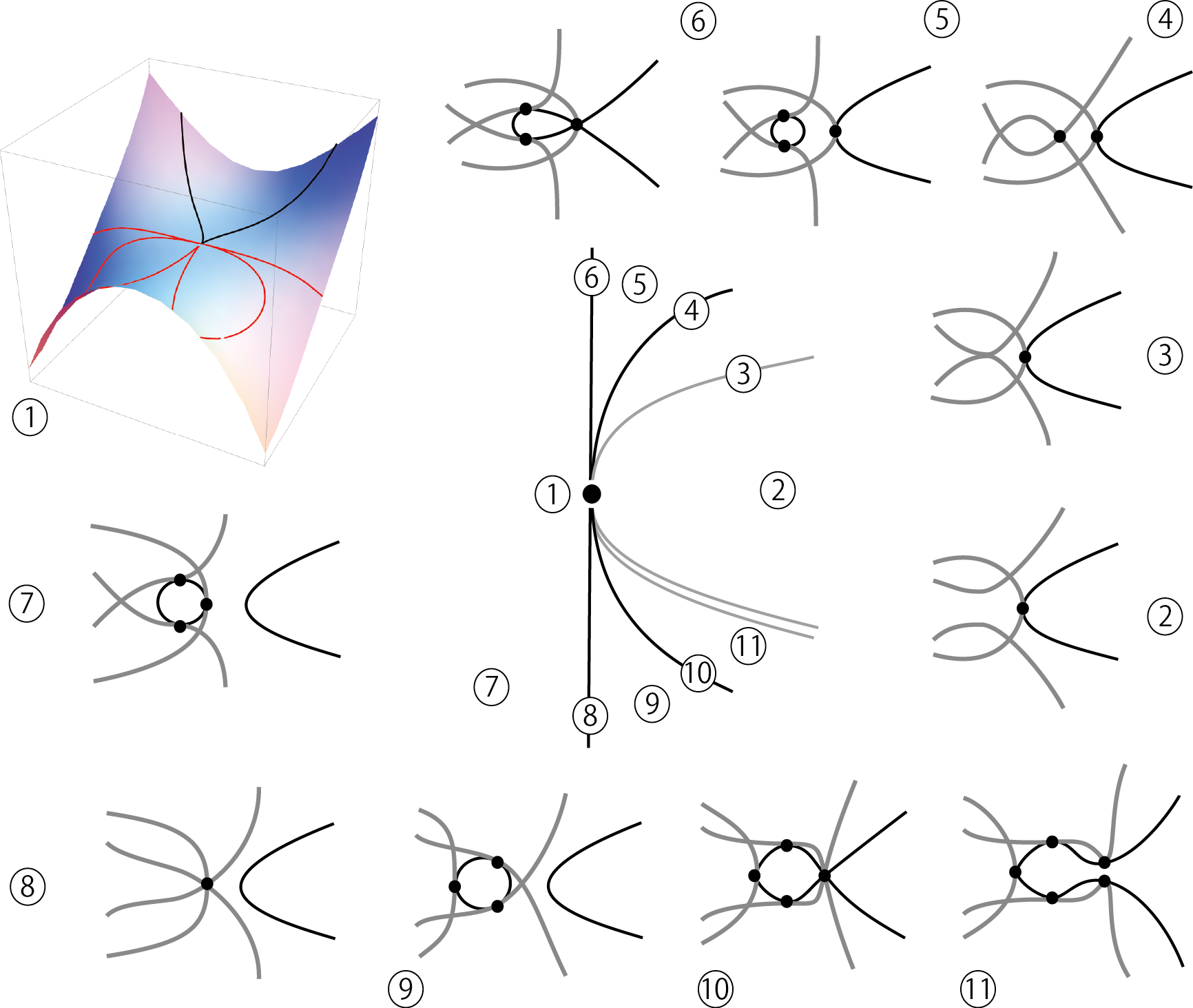}
  \caption{Bifurcation of $\Pi_2^f(-)$.}
  \label{Pif2-}
\end{figure}

%%%%%%%%%%%%%%%%%%%%%%%%%%%%%%%%%%%%%%%%%%%

\end{document}